\documentclass[12pt]{article}

\usepackage{amsmath,amssymb,amsthm,amscd,a4wide,upref}

\usepackage[enableskew]{youngtab}
\usepackage{ytableau}

\usepackage{hyperref}
\hypersetup{colorlinks,citecolor=blue,filecolor=black,linkcolor=blue,urlcolor=blue}

\usepackage{lmodern}     
\usepackage[T1]{fontenc}

\begin{document}


\newcommand{\non}{\nonumber}
\newcommand{\scl}{\scriptstyle}
\newcommand{\sclnearrow}{{\scl\nearrow}\ts}
\newcommand{\scloplus}{{\scl\bigoplus}}
\newcommand{\wt}{\widetilde}
\newcommand{\wh}{\widehat}
\newcommand{\ot}{\otimes}
\newcommand{\fand}{\quad\text{and}\quad}
\newcommand{\Fand}{\qquad\text{and}\qquad}
\newcommand{\ts}{\,}
\newcommand{\tss}{\hspace{1pt}}
\newcommand{\lan}{\langle\ts}
\newcommand{\ran}{\ts\rangle}
\newcommand{\vl}{\tss|\tss}
\newcommand{\qin}{q^{-1}}
\newcommand{\tpr}{t^{\tss\prime}}
\newcommand{\spr}{s^{\tss\prime}}
\newcommand{\di}{\partial}
\newcommand{\hra}{\hookrightarrow}
\newcommand{\antiddots}
    {\underset{\displaystyle\cdot\quad\ }
    {\overset{\displaystyle\quad\ \cdot}{\cdot}}}
\newcommand{\dddots}
    {\underset{\displaystyle\quad\ \cdot}
    {\overset{\displaystyle\cdot\quad\ }{\cdot}}}
\newcommand{\atopn}[2]{\genfrac{}{}{0pt}{}{#1}{#2}}

\newcommand{\su}{s^{}}
\newcommand{\vac}{\mathbf{1}}
\newcommand{\vacf}{\tss|0\rangle}
\newcommand{\BL}{ {\overline L}}
\newcommand{\BE}{ {\overline E}}
\newcommand{\BP}{ {\overline P}}
\newcommand{\ol}{\overline}
\newcommand{\pr}{^{\tss\prime}}
\newcommand{\ba}{\bar{a}}
\newcommand{\bb}{\bar{b}}
\newcommand{\eb}{\bar{e}}
\newcommand{\bi}{\bar{\imath}}
\newcommand{\bj}{\bar{\jmath}}
\newcommand{\bk}{\bar{k}}
\newcommand{\bl}{\bar{l}}
\newcommand{\hb}{\mathbf{h}}
\newcommand{\gb}{\mathbf{g}}
\newcommand{\For}{\qquad\text{or}\qquad}
\newcommand{\OR}{\qquad\text{or}\qquad}
\newcommand{\emp}{\mbox{}}


\newcommand{\U}{{\rm U}}
\newcommand{\Z}{{\rm Z}}
\newcommand{\ZY}{{\rm ZY}}
\newcommand{\Ar}{{\rm A}}
\newcommand{\Br}{{\rm B}}
\newcommand{\Cr}{{\rm C}}
\newcommand{\Fr}{{\rm F}}
\newcommand{\Mr}{{\rm M}}
\newcommand{\Sr}{{\rm S}}
\newcommand{\Prm}{{\rm P}}
\newcommand{\Lr}{{\rm L}}
\newcommand{\Ir}{{\rm I}}
\newcommand{\Jr}{{\rm J}}
\newcommand{\Qr}{{\rm Q}}
\newcommand{\Rr}{{\rm R}}
\newcommand{\X}{{\rm X}}
\newcommand{\Y}{{\rm Y}}
\newcommand{\DY}{ {\rm DY}}
\newcommand{\Or}{{\rm O}}
\newcommand{\SO}{{\rm SO}}
\newcommand{\GL}{{\rm GL}}
\newcommand{\Spr}{{\rm Sp}}
\newcommand{\Zr}{{\rm Z}}
\newcommand{\ev}{{\rm ev}}
\newcommand{\Pf}{{\rm Pf}}
\newcommand{\Ann}{{\rm{Ann}\ts}}
\newcommand{\Norm}{{\rm Norm\tss}}
\newcommand{\Ad}{{\rm Ad}}
\newcommand{\SY}{{\rm SY}}
\newcommand{\Pff}{{\rm Pf}\tss}
\newcommand{\Hf}{{\rm Hf}\tss}
\newcommand{\trts}{{\rm tr}\ts}
\newcommand{\otr}{{\rm otr}}
\newcommand{\row}{{\rm row}}
\newcommand{\End}{{\rm{End}\ts}}
\newcommand{\Mat}{{\rm{Mat}}}
\newcommand{\Hom}{{\rm{Hom}}}
\newcommand{\id}{{\rm id}}
\newcommand{\middd}{{\rm mid}}
\newcommand{\ch}{{\rm{ch}\ts}}
\newcommand{\ind}{{\rm{ind}\ts}}
\newcommand{\Normts}{{\rm{Norm}\ts}}
\newcommand{\mult}{{\rm{mult}}}
\newcommand{\per}{{\rm per}\ts}
\newcommand{\sgn}{{\rm sgn}\ts}
\newcommand{\sign}{{\rm sign}\ts}
\newcommand{\qdet}{{\rm qdet}\ts}
\newcommand{\sdet}{{\rm sdet}\ts}
\newcommand{\Ber}{{\rm Ber}\ts}
\newcommand{\inv}{{\rm inv}\ts}
\newcommand{\inva}{{\rm inv}}
\newcommand{\grts}{{\rm gr}\ts}
\newcommand{\grpr}{{\rm gr}^{\tss\prime}\ts}
\newcommand{\degpr}{{\rm deg}^{\tss\prime}\tss}
\newcommand{\Cond}{ {\rm Cond}\tss}
\newcommand{\Fun}{{\rm{Fun}\ts}}
\newcommand{\Rep}{{\rm{Rep}\ts}}
\newcommand{\sh}{{\rm{sh}}}
\newcommand{\weight}{{\rm{wt}\ts}}
\newcommand{\chara}{{\rm{char}\ts}}
\newcommand{\diag}{ {\rm diag}}
\newcommand{\Bos}{ {\rm Bos}}
\newcommand{\Ferm}{ {\rm Ferm}}
\newcommand{\cdet}{ {\rm cdet}}
\newcommand{\rdet}{ {\rm rdet}}
\newcommand{\imm}{ {\rm imm}}
\newcommand{\ad}{ {\rm ad}}
\newcommand{\tr}{ {\rm tr}}
\newcommand{\gr}{ {\rm gr}\tss}
\newcommand{\str}{ {\rm str}}
\newcommand{\loc}{{\rm loc}}
\newcommand{\Gr}{{\rm G}}

\newcommand{\twobar}{{\bar 2}}
\newcommand{\threebar}{{\bar 3}}


\newcommand{\AAb}{\mathbb{A}\tss}
\newcommand{\CC}{\mathbb{C}}
\newcommand{\KK}{\mathbb{K}\tss}
\newcommand{\QQ}{\mathbb{Q}\tss}
\newcommand{\SSb}{\mathbb{S}\tss}
\newcommand{\TT}{\mathbb{T}\tss}
\newcommand{\ZZ}{\mathbb{Z}\tss}
\newcommand{\Sbb}{\mathbb{S}}
\newcommand{\ZZb}{\mathbb{Z}}


\newcommand{\Ac}{{\mathcal A}}
\newcommand{\Bc}{{\mathcal B}}
\newcommand{\Cc}{{\mathcal C}}
\newcommand{\Cl}{{\mathcal Cl}}
\newcommand{\Dc}{{\mathcal D}}
\newcommand{\Ec}{{\mathcal E}}
\newcommand{\Fc}{{\mathcal F}}
\newcommand{\Jc}{{\mathcal J}}
\newcommand{\Gc}{{\mathcal G}}
\newcommand{\Hc}{{\mathcal H}}
\newcommand{\Lc}{{\mathcal L}}
\newcommand{\Nc}{{\mathcal N}}
\newcommand{\Xc}{{\mathcal X}}
\newcommand{\Yc}{{\mathcal Y}}
\newcommand{\Oc}{{\mathcal O}}
\newcommand{\Pc}{{\mathcal P}}
\newcommand{\Qc}{{\mathcal Q}}
\newcommand{\Rc}{{\mathcal R}}
\newcommand{\Sc}{{\mathcal S}}
\newcommand{\Tc}{{\mathcal T}}
\newcommand{\Uc}{{\mathcal U}}
\newcommand{\Vc}{{\mathcal V}}
\newcommand{\Wc}{{\mathcal W}}
\newcommand{\Zc}{{\mathcal Z}}
\newcommand{\HC}{{\mathcal HC}}


\newcommand{\asf}{\mathsf a}
\newcommand{\bsf}{\mathsf b}
\newcommand{\csf}{\mathsf c}
\newcommand{\nsf}{\mathsf n}


\newcommand{\Sym}{\mathfrak S}
\newcommand{\h}{\mathfrak h}
\newcommand{\q}{\mathfrak q}
\newcommand{\n}{\mathfrak n}
\newcommand{\m}{\mathfrak m}
\newcommand{\p}{\mathfrak p}
\newcommand{\gl}{\mathfrak{gl}}
\newcommand{\oa}{\mathfrak{o}}
\newcommand{\spa}{\mathfrak{sp}}
\newcommand{\osp}{\mathfrak{osp}}
\newcommand{\g}{\mathfrak{g}}
\newcommand{\kgot}{\mathfrak{k}}
\newcommand{\agot}{\mathfrak{a}}
\newcommand{\bgot}{\mathfrak{b}}
\newcommand{\sll}{\mathfrak{sl}}
\newcommand{\f}{\mathfrak{f}}
\newcommand{\z}{\mathfrak{z}}
\newcommand{\Zgot}{\mathfrak{Z}}


\newcommand{\al}{\alpha}
\newcommand{\be}{\beta}
\newcommand{\ga}{\gamma}
\newcommand{\de}{\delta}
\newcommand{\De}{\Delta}
\newcommand{\Ga}{\Gamma}
\newcommand{\ep}{\epsilon}
\newcommand{\ee}{\epsilon^{}}
\newcommand{\ve}{\varepsilon}
\newcommand{\ls}{\ts\lambda\ts}
\newcommand{\vk}{\varkappa}
\newcommand{\vs}{\varsigma}
\newcommand{\vt}{\vartheta}
\newcommand{\ka}{\kappa}
\newcommand{\vp}{\varphi}
\newcommand{\la}{\lambda}
\newcommand{\La}{\Lambda}
\newcommand{\si}{\sigma}
\newcommand{\ze}{\zeta}
\newcommand{\om}{\omega}
\newcommand{\Om}{\Omega}
\newcommand{\up}{\upsilon}


\newtheorem{thm}{Theorem}[section]
\newtheorem{lemma}[thm]{Lemma}
\newtheorem{prop}[thm]{Proposition}
\newtheorem{cor}[thm]{Corollary}
\newtheorem{conj}[thm]{Conjecture}

\theoremstyle{definition}
\newtheorem{definition}[thm]{Definition}
\newtheorem{example}[thm]{Example}

\theoremstyle{remark}
\newtheorem{remark}[thm]{Remark}

\newcommand{\bth}{\begin{thm}}
\renewcommand{\eth}{\end{thm}}
\newcommand{\bpr}{\begin{prop}}
\newcommand{\epr}{\end{prop}}
\newcommand{\ble}{\begin{lemma}}
\newcommand{\ele}{\end{lemma}}
\newcommand{\bco}{\begin{cor}}
\newcommand{\eco}{\end{cor}}
\newcommand{\bex}{\begin{example}}
\newcommand{\eex}{\end{example}}
\newcommand{\bde}{\begin{definition}}
\newcommand{\ede}{\end{definition}}
\newcommand{\bre}{\begin{remark}}
\newcommand{\ere}{\end{remark}}
\newcommand{\bcj}{\begin{conj}}
\newcommand{\ecj}{\end{conj}}

\renewcommand{\theequation}{\arabic{section}.\arabic{equation}}

\numberwithin{equation}{section}


\newcommand{\bpf}{\begin{proof}}
\newcommand{\epf}{\end{proof}}


\def\beql#1{\begin{equation}\label{#1}}

\newcommand{\bal}{\begin{aligned}}
\newcommand{\eal}{\end{aligned}}
\newcommand{\beq}{\begin{equation}}
\newcommand{\eeq}{\end{equation}}
\newcommand{\ben}{\begin{equation*}}
\newcommand{\een}{\end{equation*}}

\title{\Large\bf On the Jucys--Murphy method and fusion procedure\\ for the Sergeev superalgebra}

\author{{Iryna Kashuba,\quad Alexander Molev\quad and\quad Vera Serganova}}

\date{} 
\maketitle


\begin{abstract}
We use the Jucys--Murphy elements to
construct a complete set of primitive idempotents for the Sergeev superalgebra $\Sc_n$.
We produce seminormal forms for the simple modules over
$\Sc_n$ and over the spin symmetric group algebra
with explicit constructions of basis vectors.
We show that the idempotents can also be obtained from a new version
of the fusion procedure.


\end{abstract}

\section{Introduction}
\label{sec:int}

The {\em Sergeev superalgebra} $\Sc_n$ is the graded tensor product of two superalgebras
\beql{ser}
\Sc_n=\CC\Sym_n^-\ot \Cl_n,
\eeq
where $\CC\Sym_n^-$ is the {\em spin symmetric group algebra} generated by odd
elements $t_1,\dots,t_{n-1}$ subject to the relations
\ben
t_a^2=1,\qquad t_at_{a+1}t_a=t_{a+1}t_a t_{a+1},\qquad t_at_b=-t_bt_a,\qquad |a-b|>1,
\een
while
$\Cl_n$ is the {\em Clifford superalgebra} generated by odd elements $c_1,\dots,c_n$
subject to the relations
\ben
c_a^2=-1,\qquad c_ac_b=-c_bc_a,\qquad a\ne b.
\een
The group algebra $\CC\Sym_n$ of the symmetric group $\Sym_n$  is embedded in $\Sc_n$ so that
the adjacent transpositions $s_a=(a,a+1)\in\Sym_n$ are identified with the elements of $\Sc_n$ by
\beql{embsy}
s_a=\frac{1}{\sqrt 2}\ts t_a(c_{a+1}-c_a).
\eeq
This leads to the alternative presentation of $\Sc_n$ as the semidirect product
$\CC\Sym_n\ltimes \Cl_n$ with the relations between elements of the symmetric group
and the Clifford superalgebra given by
\beql{symgr}
s_ac_a=c_{a+1}s_a,\qquad s_ac_{a+1}=c_a s_a,\qquad s_ac_b=c_bs_a,\qquad b\ne a,a+1.
\eeq
The superalgebra $\Sc_n$ was introduced by Sergeev in \cite{s:ta}
to establish super-versions of the Schur--Weyl duality.
The
representations of $\Sc_n$ were studied in the foundational
work of Nazarov~\cite{n:ys} along with those of the {\em degenerate affine Sergeev algebra}
(also known as the {\em degenerate affine Hecke--Clifford algebra})
introduced therein. The superalgebra $\Sc_n$ is known to be semisimple
and its simple modules are parameterized by strict partitions of $n$; see e.g.
\cite[Ch.~3]{cw:dr},
\cite[Ch.~13]{k:lp} and \cite{ww:ls}
for reviews of the representation theory of $\Sc_n$.

The {\em Jucys--Murphy elements} $x_1,\dots,x_n$ of $\Sc_n$ (see \eqref{jmse} below),
introduced in \cite{n:ys} are
instrumental in the theory; they generate a commutative subalgebra of $\Sc_n$
and act semisimply on each simple module. We use these properties to
construct a complete set of primitive idempotents for $\Sc_n$; see Theorem~\ref{thm:primed}.
We apply this construction to get a new basis of $\Sc_n$ (Theorem~\ref{thm:basalg})
and explicit realizations of all its simple modules (Corollary~\ref{cor:baact}).
We also construct
all simple modules over the spin symmetric group algebra $\CC\Sym_n^-$ (Theorem~\ref{thm:spin}).
Furthermore, we give
a new version of the {\em fusion procedure} for $\Sc_n$
which yields the same primitive idempotents
by evaluating a universal rational function with values in
$\Sc_n$ (Theorem~\ref{thm:fusion}).

To explain our results in more detail, draw an analogy with
the corresponding constructions for the symmetric group which we will
now review briefly.
Recall that
the Jucys--Murphy elements $x_1,\dots,x_n$ for the symmetric group $\Sym_n$ are the sums
of transpositions
\beql{jm}
x_1=0\Fand x_a=(1,a)+\dots+(a-1, a)\qquad\text{for}\quad a=2,\dots,n
\eeq
in the group algebra $\CC\Sym_n$. They were used independently by Jucys~\cite{j:fy}
and Murphy~\cite{m:nc} to construct
the primitive idempotents $e^{}_{\Uc}$ in $\CC\Sym_n$
parameterized by the standard tableaux $\Uc$ of shapes $\la$ running over Young diagrams $\la$
with $n$ boxes. If
$\Uc$ is a standard $\la$-tableau, denote by
$\Vc$ the standard tableau
obtained from $\Uc$ by removing the box $\al$ occupied by $n$.
Then the shape of $\Vc$ is a Young diagram which we denote by $\mu$.
We let $c$ denote the content $j-i$ of the box $\al=(i,j)$.
We have the inductive formula due to \cite{j:fy}
and \cite{m:nc}:
\beql{murphyfo}
e^{}_\Uc=e^{}_\Vc\ts \frac{(x_n-a_1)\dots
(x_n-a_l)}{(c-a_1)\dots (c-a_l)},
\eeq
where $a_1,\dots,a_l$ are the contents of all addable boxes of $\mu$
except for $\al$, and
$e^{}_{\ts\tiny\young(1)}=1$ for the one-box tableau.
Note the properties
\beql{xiet}
x_a\ts e^{}_\Uc=e^{}_\Uc\ts x_a=c_a(\Uc)\ts e^{}_\Uc,
\qquad a=1,\dots,n,
\eeq
where $c_a(\Uc)$ is the content of the box in $\Uc$ occupied by $a$.

The primitive idempotents are pairwise orthogonal and form a decomposition of the identity:
\ben
e^{}_{\Uc}\tss e^{}_{\Vc}=\de^{}_{\Uc\Vc}\ts e^{}_{\Vc},\qquad
1=\sum_{\la\vdash n}\ts\sum_{\sh(\Uc)=\la} e^{}_{\Uc}.
\een
This implies the direct sum decomposition of the group algebra
\beql{decosn}
\CC\Sym_n=\bigoplus_{\la\vdash n}\ts\bigoplus_{\sh(\Uc)=\la}\ts \CC\Sym_n\tss e^{}_{\Uc}.
\eeq

Our bases for the Sergeev superalgebra will be analogous to
the seminormal bases for the group algebra $\CC\Sym_n$
constructed by Murphy~\cite[Sec.~5,6]{m:rt}. To recall them,
follow Cherednik~\cite{c:sb} to
introduce the intertwining elements $\phi_a\in\CC\Sym_n$ for $a=1,\dots,n-1$ by
\beql{phiadefsym}
\phi_a=s_a(x_a-x_{a+1})+1.
\eeq
They satisfy the braid relations (in the same form as \eqref{braidint} below)
and have the properties
\beql{phiasqasym}
\phi_a^2=1-(x_a-x_{a+1})^2.
\eeq
For each $w\in\Sym_n$ there are well-defined elements $\phi_w$ and $\phi^*_w=\phi_{w^{-1}}$
given by formulas \eqref{phiw} below. These elements were not explicitly used in \cite{m:rt},
and their properties can streamline the construction in {\em loc. cit.} Namely,
given a Young diagram $\la$ with $n$ boxes, let $\Rc^{\la}$ denote the {\em row-tableau}
of shape $\la$ obtained by writing the numbers $1,\dots,n$ consecutively
into the boxes of $\la$ by rows
from left to right, starting from the top row. For any standard $\la$-tableau $\Tc$,
there is a unique permutation $d(\Tc)\in\Sym_n$ such that $\Tc=d(\Tc)\tss \Rc^{\la}$.
For any two standard $\la$-tableaux $\Tc$ and $\Uc$
set
\beql{zetu}
\ze^{}_{\Tc\Uc}=\phi^{}_{d(\Tc)}\ts e^{}_{\Rc^{\la}}\ts \phi^*_{d(\Uc)}.
\eeq
As $\la$ runs over the Young diagrams with $n$ boxes,
the elements $\ze^{}_{\Tc\Uc}$ associated with standard
$\la$-tableaux $\Tc$ and $\Uc$
form a basis of the group algebra $\CC\Sym_n$, thus
materializing the
well-known identity
\beql{idensym}
n!=\sum_{\la\vdash n} \ts f_{\la}^2,
\eeq
where $f_{\la}$ is the number of standard tableaux of shape $\la$.
This basis was used in \cite[Thm~6.4]{m:rt} to produce {\em Young's seminormal
form} of the simple $\Sym_n$-modules. Such a module is isomorphic
to the left ideal $\CC\Sym_n\tss e^{}_{\Uc} $ occurring
in the decomposition \eqref{decosn}, where $\Uc$ is any standard $\la$-tableau.
By choosing
$\Uc$ to be the row-tableau $\Rc^{\la}$, we find that the left ideal
$\CC\Sym_n\tss e^{}_{\Rc^{\la}}$
is spanned by the basis elements
\beql{zedefsym}
\ze^{}_{\Tc}:=\phi^{}_{d(\Tc)}\ts e^{}_{\Rc^{\la}}=
e^{}_{\Tc}\ts\phi^{}_{d(\Tc)},
\eeq
parameterized by the standard $\la$-tableaux $\Tc$.
The properties \eqref{xiet} imply that $x_a\tss \ze^{}_{\Tc}=c_a(\Tc)\tss \ze^{}_{\Tc}$.
If the tableau $s_a\Tc$ obtained from $\Tc$ by swapping $a$ and $a+1$ is standard and
the length of the permutation $d(s_a\Tc)$ is greater than
that of $d(\Tc)$, then $\phi_a\tss \ze^{}_{\Tc}=\ze^{}_{s_a\Tc}$.
By using \eqref{phiadefsym}
together with \eqref{phiasqasym}, we get
explicit formulas for the action of the generators $s_a$ on the
basis vectors $\ze^{}_{\Tc}$.

Returning to the primitive idempotents $e^{}_\Uc$, note that
another way to express them
is provided by the {\em fusion procedure}
originated from Jucys~\cite{j:yo}.
It was
re-discovered by Cherednik~\cite{c:sb}; detailed
arguments were also given by Nazarov~\cite{n:yc}.
Take $n$ variables $u_1,\dots,u_n$
and consider the rational function with values in $\CC\Sym_n$
defined by
\beql{phiuuusym}
\Phi(u_1,\dots,u_n)=
\prod_{1\leqslant a<b\leqslant n}\Big(1-\frac{(a,b)}{u_a-u_b}\Big),
\eeq
where the product is taken in
the lexicographical order on the set of pairs $(a,b)$.
Set $c_a=c_a(\Uc)$ for $a=1,\dots,n$. According to the version given in
\cite{m:fp}, the consecutive evaluations
of the rational function \eqref{phiuuusym} at the contents are well-defined
and the value coincides with the primitive idempotent $e^{}_\Uc$
up to a constant factor:
\beql{fusn}
\Phi(u_1,\dots,u_n)\big|_{u_1=c_1}
\big|_{u_2=c_2}\dots \big|_{u_n=c_n}=\frac{n!}{f_{\la}}\ts e^{}_\Uc.
\eeq

The factors in \eqref{phiuuusym} evaluated in the context of the Schur--Weyl duality
coincide with the {\em Yang $R$-matrices}.
We refer the reader e.g. to \cite[Ch.~6,7]{m:yc} for applications of the fusion formula
\eqref{fusn} to representations of the general linear Lie algebra and associated Yangian
along with a brief review of the literature on the extension of the fusion procedure to other classes
of algebras.

Our results for the Sergeev superalgebra $\Sc_n$
and the spin symmetric group algebra $\CC\Sym_n^-$
are parallel to those for
the group algebra $\CC\Sym_n$. We define elements $e^{}_\Uc$
of $\Sc_n$ by formula \eqref{murphyse} analogous to
\eqref{murphyfo}, where the role of the standard tableaux is now played by
the {\em standard shifted barred tableaux} $\Uc$ as defined in Sec.~\ref{sec:jm}. We prove that
the elements $e^{}_\Uc$ are idempotents in $\Sc_n$. We show that their
enhancements via idempotents of the Clifford superalgebra yield a complete
set of pairwise orthogonal primitive idempotents thus producing a decomposition
of $\Sc_n$ into a direct sum of left ideals \eqref{sedec} analogous to \eqref{decosn}.
The left ideals are simple $\Sc_n$-modules.
This can be regarded
as a refinement of the construction of Nazarov,
where left ideals of $\Sc_n$ were produced
with the use of a fusion procedure \cite[Thm~7.2]{n:ys}
whose version was also given in \cite{n:ci}.
They were shown to bear a direct sum
of several copies of simple modules.

As in \cite{m:rt}, we use
the primitive idempotents to construct a basis of
the Sergeev superalgebra parameterized by pairs of
standard shifted tableaux and subsets of barred boxes
by analogy with \eqref{zetu} (Theorem~\ref{thm:basalg}).
We rely on the properties of the intertwiners $\phi_a$ as defined in \eqref{phiadef} which were
originally introduced in \cite{n:ys}. They are also used in our
derivation of the seminormal form of the simple modules over $\Sc_n$
and their explicit construction as left ideals
of $\Sc_n$ (Corollary~\ref{cor:baact}).

By the results of Brundan and Kleshchev~\cite[Thm~3.4]{bk:pr},
the representations of the superalgebras $\Sc_n$ and $\CC\Sym_n^-$ are closely
related to each other by an adjoint pair of exact functors.
Using this relationship, we
restrict certain representations of $\Sc_n$ to the spin symmetric group
algebra $\CC\Sym_n^-$ to produce
the seminormal form for all simple $\CC\Sym_n^-$-modules $V^{\la}$
(Theorem~\ref{thm:spin}). A different construction of these modules
was given earlier by Nazarov~\cite{n:yo}. Note also the related work of
Vershik and Sergeev~\cite{vs:na}, where
formulas for
the action of the generators of $\CC\Sym_n^-$
on certain vectors parameterized by standard $\la$-tableaux were given.
They correspond to formulas \eqref{taspin} below for the vectors $v^{}_{\Tc}$,
but the explicit basis construction requires an additional step
involving the Clifford superalgebra. To get the seminormal form of the simple module $V^{\la}$
in Theorem~\ref{thm:spin},
we use the `normalized' vectors $\theta^{}_{\Tc}$ defined in \eqref{thetadef}
and the simple right modules over the Clifford superalgebras.

Our version of the fusion procedure
for $\Sc_n$ is analogous to \eqref{fusn}. It builds
on the work of Nazarov~\cite{n:ys, n:ci} and
uses the rational function
$\Phi(u_1,\dots,u_n)$ with values in $\Sc_n$ given in \eqref{phiuuu}
which was introduced in {\em loc. cit.}
The use of the {\em consecutive evaluations}
of the variables as in \eqref{fusn} provides a new version
of the fusion procedure
similar to \cite{m:fp} thus showing its equivalence to
the Jucys--Murphy-type formulas in \eqref{murphyse} below.

We expect that
the Schur--Sergeev duality between $\Sc_n$ and the queer Lie superalgebra $\q_N$
as reviewed e.g. in
\cite[Ch.~3]{cw:dr}, should make the idempotent formulas provided by the Jucys--Murphy method
and fusion procedure useful
in the representation theory of $\q_N$ and associated objects.

The first version of our preprint appeared in the arXiv almost simultaneously with
an independent work of Shuo Li and Lei Shi~\cite{ls:sb}. They investigate
more general classes of objects, while in the particular case of the Sergeev superalgebra, 
their results concerning primitive idempotents and seminormal forms
are consistent with ours, as already pointed out in
\cite[Remark~5.9]{ls:sb}.

\section{Jucys--Murphy method}
\label{sec:jm}

Introduce analogs of the transpositions in the Sergeev superalgebra \eqref{ser}
by setting
\ben
t_{ab}=(-1)^{b-a-1}\ts t_{b-1}\dots t_{a+1}t_a t_{a+1}\dots t_{b-1},\qquad a<b,
\een
and $t_{ba}=-t_{ab}$.
The odd Jucys--Murphy elements \cite{s:hd} of $\Sc_n$ are given by
\ben
m_1=0,\qquad m_a=t_{1a}+\dots+t_{a-1,a},\qquad a=2,\dots,n,
\een
and the even Jucys--Murphy elements \cite{n:ys} are
\beql{jmse}
x_a=\sqrt{2}\ts m_a\tss c_a,\qquad a=1,\dots,n.
\eeq
Note that $x_a^2=2\tss m_a^2$ and the $x_a$ pairwise commute.

The irreducible representations of both superalgebras $\Sc_n$ and $\CC\Sym_n^-$ are
parameterized by strict partitions
$\la=(\la_1,\dots,\la_{\ell})$ of $n$ with $\la_1>\dots>\la_{\ell}>0$
and $\la_1+\dots+\la_{\ell}=n$. We will set $\ell(\la)=\ell$ to denote the length
of the partition and
write $\la\Vdash n$ to indicate that $\la$ is a strict partition on $n$.
A strict partition $\la$ is usually depicted by the shifted Young diagram obtained
from the usual Young diagram by shifting row $i$ to the right by $i-1$ unit boxes,
as illustrated by the diagram for $\la=(6,3,1)$:
\ben
\ydiagram{6, 1+3, 2+1}\ .
\een

A ({\em shifted}) $\la$-{\em tableau} is obtained by writing the numbers $1,\dots,n$ bijectively
into the boxes of the shifted Young diagram $\la$.
Such a tableau
is called {\em standard} if its entries increase
from left to right in each row and from top to bottom in each column.
The entries in the boxes $(i,i)$ are called the {\em diagonal entries}.
For instance,
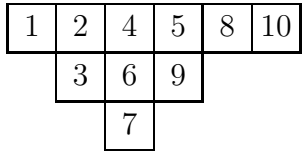
\ben
\begin{ytableau}
    1 & 2 & 4 & 5 & 8 & 10\\
    \none & 3 & 6 & 9\\
    \none & \none & 7
\end{ytableau}\phantom{\ .}
\een
is a standard $\la$-tableau for $\la=(6,3,1)$ with the diagonal entries $1,3$ and $7$.

The respective dimensions
of the simple $\Sc_n$-module $U^{\la}$ and the simple $\CC\Sym_n^-$-module $V^{\la}$ are given by
\ben
\dim U^{\la}=2^{n-\lfloor \frac{\ell(\la)}{2}\rfloor}_{}\ts g_{\la}\Fand
\dim V^{\la}=2^{\lceil \frac{n-\ell(\la)}{2}\rceil}_{}\ts g_{\la},
\een
where $g_{\la}$ is the number of standard $\la$-tableaux
which is found by the Schur formula
\beql{gla}
g_{\la}=\frac{n!}{\la_1!\dots \la_{\ell}!}\ts \prod_{1\leqslant i<j\leqslant \ell}\ts
\frac{\la_i-\la_j}{\la_i+\la_j}.
\eeq
Both superalgebras are semisimple;
the corresponding Wedderburn decomposition for $\Sc_n$ involves direct summands
of types $\mathsf{M}$ and $\mathsf{Q}$ which are
associated with the modules $U^{\la}$, where the length
$\ell(\la)$ is even and odd,
respectively:
\ben
\Sc_n\cong\bigoplus_{\la\Vdash n,\ \ell(\la)\ts  \text{even}}\ts
\Mat\big(2^{n-p}g_{\la}|2^{n-p}g_{\la}\big)
\oplus\bigoplus_{\la\Vdash n,\ \ell(\la)\ts  \text{odd}}\ts \Qr\big(2^{n-p}g_{\la}\big),
\een
where
$p=\lfloor \frac{\ell(\la)}{2}\rfloor+1$;
see e.g. \cite[Ch.~3]{cw:dr}
for proofs.
The dimension count in the decomposition
yields the Schur identity analogous to \eqref{idensym}:
\beql{idendeim}
n!=\sum_{\la\Vdash n} \ts 2^{n-\ell(\la)}\ts g_{\la}^2.
\eeq

We need to extend the set of tableaux by allowing any non-diagonal entry to
occur with a bar on it, as illustrated for the above standard tableau:
\ben
\begin{ytableau}
    1 & 2 & \ol{4} & \ol{5} & 8 & \ol{10}\\
    \none & 3 & \ol{6} & 9\\
    \none & \none & 7
\end{ytableau}
\ .
\een
Note that given
a strict
partition $\la\Vdash n$, the number of the
corresponding
standard barred tableaux equals
$2^{n-\ell(\la)}\tss g_{\la}$.

We will use symbols $\asf, \bsf, \csf$, etc.
to denote arbitrary barred or unbarred
entries of $\Uc$.
Given a standard barred tableaux $\Uc$, introduce
the {\em signed content} $\ka_{\asf}(\Uc)$
of any barred or unbarred entry $\asf$ of $\Uc$
by the formula
\ben
\ka_{\asf}(\Uc)=\begin{cases}\phantom{-\ts}\sqrt{\si_a(\si_a+1)}
\qquad&\text{if $\asf=a$ is unbarred},\\[0.4em]
-\ts\sqrt{\si_a(\si_a+1)}\qquad&\text{if $\asf=\bar a$ is barred},
\end{cases}
\een
where $\si_a=j-i$ is the content of the box $(i,j)$ of $\la$
occupied by $a$ or $\bar a$. Similarly, we will
extend the notation for the Jucys--Murphy elements
by setting $x_{\asf}=x_{a}$ for $\asf=a$ and $\asf=\bar a$.

By analogy with \eqref{murphyfo},
for any standard barred tableau $\Uc$
introduce the element $e^{}_{\Uc}$ of the Sergeev superalgebra $\Sc_n$ by induction, setting
$e^{}_{\ts\tiny\young(1)}=1$ for the one-box tableau, and
\beql{murphyse}
e^{}_{\Uc}=e^{}_{\Vc}\ts \frac{(x_n-b_1)\dots
(x_n-b_p)}{(\ka-b_1)\dots (\ka-b_p)},
\eeq
where $\Vc$ is the barred tableau
obtained from $\Uc$ by removing the box $\al$ occupied by $n$ (resp. $\bar n$).
The shape of $\Vc$ is a shifted diagram which we denote by $\mu$ and
$b_1,\dots,b_p$ are the signed contents in all addable boxes of $\mu$ (barred and unbarred),
except for the entry $n$ (resp. $\bar n$), while $\ka$ is the signed content of
the entry $n$ (resp. $\bar n$).

\bex\label{ex:smdi3}
For $n=3$ we have six barred tableaux $\Uc$:
\ben
\begin{ytableau}
    1 & 2 & 3
\end{ytableau}
\qquad
\begin{ytableau}
    1 & 2 & \ol 3
\end{ytableau}
\qquad\begin{ytableau}
    1 & \ol 2 & 3
\end{ytableau}
\qquad
\begin{ytableau}
    1 & \ol 2 & \ol 3
\end{ytableau}
\qquad
\begin{ytableau}
    1 & 2\\
    \none&3
\end{ytableau}
\qquad
\begin{ytableau}
    1 & \ol 2\\
    \none&3
\end{ytableau}\ \ .
\een
The respective elements $e^{}_{\Uc}\in\Sc_3$ are
\ben
\frac{\sqrt{2}+x_2}{2\tss\sqrt{2}}\cdot \frac{x_3(x_3+\sqrt{6})}{12},\qquad
\frac{\sqrt{2}+x_2}{2\tss\sqrt{2}}\cdot \frac{x_3(x_3-\sqrt{6})}{12},\qquad
\frac{\sqrt{2}-x_2}{2\tss\sqrt{2}}\cdot \frac{x_3(x_3+\sqrt{6})}{12},
\een
\ben
\frac{\sqrt{2}-x_2}{2\tss\sqrt{2}}\cdot \frac{x_3(x_3-\sqrt{6})}{12},\qquad
\frac{\sqrt{2}+x_2}{2\tss\sqrt{2}}\cdot \frac{6-x_3^2}{6},\qquad
\frac{\sqrt{2}-x_2}{2\tss\sqrt{2}}\cdot \frac{6-x_3^2}{6}.
\een

\noindent
Note that $x_2^2=2$ and $x_3^3=6x_3$ so that all these elements are idempotents in $\Sc_3$.
\eex

\bpr\label{prop:idemp}
All elements $e^{}_{\Uc}$ are idempotents in $\Sc_n$. They are
pairwise orthogonal and form a decomposition of the identity:
\ben
e^{}_{\Uc}\tss e^{}_{\Vc}=\de^{}_{\Uc\Vc}\ts e^{}_{\Vc},\qquad
1=\sum_{\la\Vdash n}\ts\sum_{\sh(\Uc)=\la} e^{}_{\Uc}.
\een
Moreover,
\beql{xietsign}
x_{\asf}\ts e^{}_{\Uc}=e^{}_{\Uc}\ts x_{\asf}=\ka_{\asf}(\Uc)\ts e^{}_{\Uc},
\eeq
where $\ka_{\asf}(\Uc)$ is the signed content of the entry $\asf$ of\ \ts $\Uc$.
\epr

\bpf
We will use the explicit construction of $\Sc_n$-modules going back to
\cite{n:ys} as developed in \cite{hks:da}, \cite{vs:na} and \cite{w:cs}; see also
its review in \cite{ww:ls}. For each strict partition $\la\Vdash n$
consider the module over $\Sc_n$ afforded by the vector space
\beql{ula}
\wh U^{\la}=\bigoplus_{\sh(\Tc)=\la}\ts \Cl_n \tss v^{}_{\Tc}
\eeq
with the basis vectors $v^{}_{\Tc}$ associated with the (unbarred) standard $\la$-tableaux $\Tc$.
The action of generators of $\Sc_n$ on the basis vectors is given
by explicit formulas as in \cite[Sec.~5.1]{hks:da} and \cite[Sec.~7.5]{ww:ls}:
\begin{align}
s_a\tss v^{}_{\Tc}&=\Bigg(\frac{1}{\ka_{a+1}(\Tc)-\ka_{a}(\Tc)}
+\frac{c_ac_{a+1}}{\ka_{a+1}(\Tc)+\ka_{a}(\Tc)}\Bigg)\tss v^{}_{\Tc}
+\Yc_a(\Tc)\ts v^{}_{s_a\Tc},
\label{saavt}\\[0.5em]
x_a\tss v^{}_{\Tc}&=\ka_{a}(\Tc)\ts v^{}_{\Tc},
\label{jmacti}
\end{align}
where we assume that $v^{}_{s_a\Tc}=0$ if the tableau $s_a\Tc$ is not standard, and
\beql{yatc}
\Yc_a(\Tc)=\sqrt{A\big(\ka_a(\Tc),\ka_{a+1}(\Tc)\big)},
\eeq
with
\beql{auv}
A(u,v)=1-\frac{1}{(u-v)^2}-\frac{1}{(u+v)^2}.
\eeq
The formulas for the Jucys--Murphy elements in \eqref{jmacti} follow by taking
the images of the respective generators of the degenerate affine
Hecke--Clifford algebra in $\Sc_n$; see {\em loc. cit.}
Now observe that the representation
\beql{uhat}
\wh U=\bigoplus_{\la\Vdash n}  \wh U^{\la}
\eeq
of the superalgebra $\Sc_n$ is faithful.
Indeed, as a consequence of \cite[Thm~5.13]{hks:da} or \cite[Thm~7.17]{ww:ls},
each $\Sc_n$-module $\wh U^{\la}$ is isomorphic to a direct sum
of $2^{\lfloor \ell(\la)/2\rfloor}_{}$
copies of simple modules $U^{\la}$. If an element $y\in \Sc_n$ acts as zero
in $\wh U$, then $y$ acts as zero
in every simple $\Sc_n$-module. This means that $y$ belongs to the radical of $\Sc_n$,
implying $y=0$ because the superalgebra $\Sc_n$ is semisimple.

To complete the proof,
it suffices to verify the desired properties of the elements of $\Sc_n$ for their images
in the representation $\wh U$.
Due to the relations
\beql{xaca}
x_ac_a=-c_ax_a,\qquad x_a\tss c_b=c_b\tss x_a,\qquad a\ne b,
\eeq
in $\Sc_n$, each subspace
$\Cl_n\tss v^{}_{\Tc}$ is stable under the action of the subalgebra of $\Sc_n$, generated
by $x_1,\dots,x_n$ and $c_1,\dots,c_n$. Therefore,
it suffices to consider the action of the elements of $\Sc_n$
in the statement of the proposition on
these subspaces.
Along with the basis vectors $v^{}_{\Tc}$ of $\wh U^{\la}$ consider the vectors
$v^{}_{\Uc}$ associated with the standard barred tableaux $\Uc$. Namely,
suppose that a certain tableau $\Tc$ is obtained from $\Uc$ by unbarring all
barred entries $\bar a_1,\dots,\bar a_r$ with $a_1<\dots< a_r$.
Then set
\beql{vbar}
v^{}_{\ts\Uc}=c_{a_1}\dots c_{a_r}\ts v^{}_{\Tc}.
\eeq
Observe that
\ben
x_{\asf}\ts v^{}_{\Uc}=\ka_{\asf}(\Uc)\ts v^{}_{\Uc}
\een
for any entry $\asf$ of $\Uc$.
Hence the definition \eqref{murphyse} implies that
in $\wh U$ we have
\beql{idema}
e^{}_{\Uc}\ts v^{}_{\Vc}=\de^{}_{\Uc\tss\Vc}\ts v^{}_{\Vc}
\eeq
for any two barred tableaux $\Uc$ and $\Vc$ with $n$ boxes. Indeed,
suppose that $e^{}_{\Uc}\ts v^{}_{\Vc}\ne 0$.
Since the Jucys--Murphy elements pairwise commute, an easy
induction on $n$ shows that $\Uc=\Vc$ and
$e^{}_{\Uc}\ts v^{}_{\Uc}=v^{}_{\Uc}$
thus verifying \eqref{idema}.

Furthermore, if $\Uc$ is a standard barred tableau of shape $\la$, suppose that
its diagonal entries are $d_1<\dots<d_{\ell}$.
By \eqref{murphyse} the idempotent $e^{}_{\Uc}$ depends on
the squares $x^2_{d_1},\dots,x^2_{d_{\ell}}$ of the corresponding Jucys--Murphy
elements. Then \eqref{xaca} implies that $e^{}_{\Uc}$ commutes
with the Clifford generators $c_{d_1},\dots,c_{d_{\ell}}$.
On the other hand, if $a$ is an unbarred non-diagonal entry of $\Uc$, then
\eqref{murphyse} and \eqref{xaca} imply the relations
\beql{flibar}
e^{}_{\Uc}\ts c_a=c_a\ts e^{}_{\Uc'}\Fand e^{}_{\Uc'}\ts c_a=c_a\ts e^{}_{\Uc},
\eeq
where the tableau $\Uc'$ is obtained from $\Uc$ by replacing $a$ with $\bar a$.
It follows from these observations that
$e^{2}_{\Uc}=e^{}_{\Uc}$, because both sides act in the same way
on all subspaces $\Cl_n\tss v^{}_{\Tc}$. The remaining claims follow by the same argument.
\epf

By Proposition~\ref{prop:idemp}, we have the direct sum decomposition
into left ideals
\ben
\Sc_n=\bigoplus_{\la\Vdash n}\ts\bigoplus_{\sh(\Uc)=\la}\ts \Sc_n\tss e^{}_{\Uc},
\een
with the second sum taken over standard barred tableaux $\Uc$ of shape $\la$.
It is implied by the proof that
each left ideal $\Sc_n\tss e^{}_{\Uc}$ is a direct sum of copies of the
corresponding simple modules $U^{\la}$; cf. \cite[Thm~8.3]{n:ys}.
We will refine this decomposition by introducing idempotents for the Clifford superalgebra.
Recall the isomorphisms
\beql{cla}
\Cl_{2m}\cong\Mat(2^{m-1}|2^{m-1})\Fand \Cl_{2m+1}\cong\Qr(2^m).
\eeq
Each of $\Cl_{2m}$ and $\Cl_{2m+1}$ contains $2^m$ primitive idempotents
$
\Ec_1,\dots,\Ec_{2^m}
$
corresponding
to the diagonal matrix units in $\Mat(2^{m-1}|2^{m-1})$ and $\Qr(2^m)$, respectively.
They are pairwise orthogonal and form a decomposition of the identity:
\ben
\Ec_a\ts \Ec_b=\de_{ab}\ts \Ec_b,\qquad 1=\sum_{a=1}^{2^m} \Ec_a.
\een

As we pointed out in the proof of Proposition~\ref{prop:idemp},
if a standard barred tableau $\Uc$ of shape $\la$ has
diagonal entries $d_1<\dots<d_{\ell}$, then $e^{}_{\Uc}$ commutes
with the Clifford generators $c_{d_1},\dots,c_{d_{\ell}}$.
Consider the subalgebra $\Cl^{\Uc}_{\ell}$ of $\Cl_n$ generated by $c_{d_1},\dots,c_{d_{\ell}}$,
which is isomorphic to the Clifford superalgebra $\Cl_{\ell}$. Introduce the corresponding
idempotents $\Ec^{\Uc}_1,\dots,\Ec^{\Uc}_{2^m}$ in $\Cl^{\Uc}_{\ell}$,
where $m$ is defined by $\ell=2m$ or $\ell=2m+1$ for the even and odd $\ell$,
respectively. We can conclude that the products
\ben
e^{(r)}_{\Uc}:=\Ec^{\Uc}_r\tss e^{}_{\Uc}= e^{}_{\Uc}\tss\Ec^{\Uc}_r,
\een
with $r=1,\dots,2^m$,
are idempotents in $\Sc_n$.

\bth\label{thm:primed}
We have
the direct sum decomposition
\beql{sedec}
\Sc_n=\bigoplus_{\la\Vdash n}\ts\bigoplus_{\sh(\Uc)=\la}\ts
\bigoplus_{r=1}^{\ 2^m}\ts \Sc_n\tss e^{(r)}_{\Uc},
\eeq
with $m=\lfloor \frac{\ell(\la)}{2}\rfloor$.
Moreover, the left ideal $\Sc_n\ts e^{(r)}_{\Uc}$ is isomorphic to
the irreducible $\Sc_n$-module $U^{\la}$.
\eth

\bpf
The elements $e^{(r)}_{\Uc}$ with $\Uc$ running
over standard barred tableaux $\Uc$ with $n$ boxes and
$r=1,\dots,2^m$,
where $\la$ is the shape of $\Uc$, are pairwise orthogonal
idempotents in $\Sc_n$. It follows from Proposition~\ref{prop:idemp}
that they form a decomposition of the identity which implies the first part
of the theorem.

Furthermore, it is clear from \eqref{idema} that each left ideal $\Sc_n\ts e^{(r)}_{\Uc}$
is nonzero and isomorphic to a direct sum of some copies of the irreducible $\Sc_n$-module $U^{\la}$.
Hence,
\beql{indecdi}
\dim \Sc_n\tss e^{(r)}_{\Uc}\geqslant 2^{n-\lfloor\frac{\ell(\la)}{2}\rfloor}\ts g_{\la}.
\eeq
Now use the decomposition \eqref{sedec}.
The number of standard barred $\la$-tableaux equals $2^{n-\ell(\la)}\ts g_{\la}$
so that the dimension of the module on the right hand side
is at least
\ben
\sum_{\la\Vdash n}\ts 2^{n-\ell(\la)}\ts g_{\la}
\times 2^{\lfloor\frac{\ell(\la)}{2}\rfloor}
\times 2^{n-\lfloor\frac{\ell(\la)}{2}\rfloor}\ts g_{\la}
=\sum_{\la\Vdash n} \ts 2^{2\tss n-\ell(\la)}\ts g_{\la}^2.
\een
However, this coincides with $\dim \Sc_n=2^n\ts n!$ by the identity in \eqref{idendeim}.
Therefore, the inequalities in \eqref{indecdi} are, in fact, equalities.
This proves that the left ideals are irreducible modules.
\epf

It is clear from the proof of the theorem that
for any standard barred tableaux $\Uc$ of shape $\la$
the left ideal
$\Sc_n\tss e^{}_{\Uc}$ is isomorphic to the direct sum of
$2^{\lfloor \ell(\la)/2\rfloor}_{}$
copies of simple modules $U^{\la}$. Hence, we have an isomorphism of
$\Sc_n$-modules
\beql{isomsn}
\wh U^{\la}\cong\Sc_n\tss e^{}_{\Uc}.
\eeq

Note also that, as an obvious consequence of Proposition~\ref{prop:idemp}, we get the expression
for the Jucys--Murphy element $x_n$
in terms of the idempotents $e^{}_{\Uc}$:
\beql{xnid}
x_n=\sum_{\la\Vdash n}\ts\sum_{\sh(\Uc)=\la} \ka_{\nsf}(\Uc)\ts e^{}_{\Uc},
\eeq
where $\ka_{\nsf}(\Uc)$ is the signed content of the entry $\nsf=n$ or $\nsf=\bar n$ of $\Uc$.

\section{Bases in simple modules}
\label{sec:sf}

We will produce explicit constructions
of the simple modules over the Sergeev superalgebra $\Sc_n$ and the spin symmetric
group algebra $\CC\Sym_n^-$; as before, we regard the latter as a superalgebra.

\subsection{Seminormal form for the Sergeev superalgebra}
\label{subsec:ser}

Following \cite{n:ys}, for $a=1,\dots,n-1$
introduce the intertwining elements $\phi_a$ of $\Sc_n$ by
\beql{phiadef}
\phi_a=s_a(x_a^2-x_{a+1}^2)+x_a+x_{a+1}-c_ac_{a+1}(x_a-x_{a+1}).
\eeq
They satisfy the braid relations
\beql{braidint}
\phi_a\phi_{a+1}\phi_a=\phi_{a+1}\phi_a \phi_{a+1},\qquad \phi_a\phi_b=\phi_b\phi_a,\qquad |a-b|>1,
\eeq
and have the properties
\beql{phiasqa}
\phi_a^2=2(x_a^2+x_{a+1}^2)-(x_a^2-x_{a+1}^2)^2
\eeq
together with
\beql{phia}
\phi_a c_a=c_{a+1}\phi_a,\qquad \phi_a c_{a+1}=c_a\phi_a,\qquad \phi_a c_b=c_b\phi_a,\quad b\ne a,a+1.
\eeq

We will need the formulas for the action of the intertwiners $\phi_a$
on the basis vectors $v^{}_{\Tc}$ of the $\Sc_n$-module $\wh U^{\la}$
introduced in \eqref{ula}. They
are implied by \eqref{saavt} and \eqref{jmacti} and
have the form
\beql{phiaa}
\phi_a\tss v^{}_{\Tc}=\big(\ka_a(\Tc)^2-\ka_{a+1}(\Tc)^2\big)\Yc_a(\Tc)\tss v^{}_{s_a\Tc},
\eeq
where $\Yc_a(\Tc)$ is defined in \eqref{yatc} and
we assume that $v^{}_{s_a\Tc}=0$ if the tableau $s_a\Tc$
obtained from $\Tc$ by swapping $a$ and $a+1$ is not standard.

For each $w\in\Sym_n$ there are well-defined elements $\phi_w$ and $\phi^*_w=\phi_{w^{-1}}$
given by
\beql{phiw}
\phi_w=\phi_{a_1}\dots \phi_{a_r}\Fand \phi^*_w=\phi_{a_r}\dots \phi_{a_1},
\eeq
where $w=s_{a_1}\dots s_{a_r}$ is a reduced decomposition.
Relations \eqref{phia} imply that
\beql{phicli}
\phi_w c_a=c_{w(a)} \phi_w.
\eeq

For any $\la\Vdash n$,
the symmetric group $\Sym_n$ acts naturally on the set of barred $\la$-tableaux.
Namely, if $a$ (resp. $\bar a$) is the entry of a tableau
$\Uc$ in the box $\al$ and $w\in\Sym_n$, then the entry of the tableau $w\tss\Uc$
in the box $\al$ is $w(a)$ (resp. $\overline{w(a)}$).

\ble\label{lem:inide}
For any standard barred $\la$-tableau $\Uc$ and $a\in\{1,\dots,n-1\}$ we have
\beql{phiau}
\phi_a\tss e^{}_{\Uc}=e^{}_{s_a\Uc}\tss\phi_a,
\eeq
where we assume that $e^{}_{s_a\tss\Uc}=0$ if the tableau $s_a\tss\Uc$ is not standard.
\ele

\bpf
Suppose first that all entries of $\Uc$ are unbarred.
As we pointed out in the proof of Proposition~\ref{prop:idemp},
the representation $\wh U$ of $\Sc_n$ is faithful. Therefore,
it suffices to verify that the relation holds for the images of the elements on
both sides of \eqref{phiau} acting in $\wh U$.
However, due to \eqref{idema} and \eqref{phiaa}, both sides of \eqref{phiau}
act in the same way on any basis vector $v^{}_{\Tc}$.

The extension of the claim to arbitrary standard barred tableaux $\Uc$
is straightforward from \eqref{flibar} and \eqref{phia}.
\epf

For every shifted Young diagram $\la$ with $n$ boxes we let $\Rc^{\la}$ denote the {\em row-tableau}
of shape $\la$ obtained by writing the numbers $1,\dots,n$ consecutively
into the boxes of $\la$ by rows
from left to right, starting from the top row. Given any standard $\la$-tableau $\Tc$,
there is a unique permutation $d(\Tc)\in\Sym_n$ such that $\Tc=d(\Tc)\tss \Rc^{\la}$.

Furthermore, let $\be$ be a subset of the set of non-diagonal boxes of $\la$.
Denote by $\Rc^{\la,\be}$ the tableau obtained from $\Rc^{\la}$ by adding the bars to all entries
occupying the boxes of $\be$. For any two standard unbarred tableaux $\Tc$ and $\Uc$
of shape $\la$ set
\ben
\ze^{\be}_{\Tc\Uc}=\phi^{}_{d(\Tc)}\ts e^{}_{\Rc^{\la,\be}}\ts \phi^*_{d(\Uc)}.
\een
Lemma~\ref{lem:inide} implies the following equivalent expressions:
\beql{zeequi}
\ze^{\be}_{\Tc\Uc}=e^{}_{\ol\Tc}\ts \phi^{}_{d(\Tc)}\ts {\phi^*_{d(\Uc)}}=
\phi^{}_{d(\Tc)}\ts \phi^*_{d(\Uc)}\ts e^{}_{\ol\Uc},
\eeq
where $\ol\Tc$ and $\ol\Uc$ are the tableaux obtained from
$\Tc$ and $\Uc$, respectively, by adding the bars to all entries
occupying the boxes of $\be$. Since $e^{2}_{\ol\Uc}=e^{}_{\ol\Uc}$, we can also write
\beql{zeeqet}
\ze^{\be}_{\Tc\Uc}=e^{}_{\ol\Tc}\ts\phi^{}_{d(\Tc)}\ts \phi^*_{d(\Uc)}\ts e^{}_{\ol\Uc}
=\phi^{}_{d(\Tc)}\ts e^{}_{\Rc^{\la,\be}}\ts \phi^*_{d(\Uc)}\ts e^{}_{\ol\Uc}.
\eeq

The next theorem shows that the elements $\ze^{\be}_{\Tc\Uc}$ can be regarded
as Sergeev superalgebra analogues of the seminormal basis vectors constructed
by Murphy~\cite[Sec.~5,6]{m:rt} for the Hecke algebras.

\bth\label{thm:basalg}
As $\la$ runs over shifted Young diagrams with $n$ boxes,
the elements $\ze^{\be}_{\Tc\Uc}$ associated with standard
unbarred tableaux $\Tc$ and $\Uc$
of shape $\la$ and sets $\be$ of non-diagonal boxes of $\la$ form a
basis of the Sergeev superalgebra $\Sc_n$ over $\Cl_n$.
\eth

\bpf
The number of elements is
\ben
\sum_{\la\Vdash n}\ts 2^{n-\ell(\la)}\ts  g_{\la}^2,
\een
which equals $n!$ by \eqref{idendeim} and so coincides with the rank of $\Sc_n$ as a $\Cl_n$-module.
Therefore, it is enough to show that the elements are linearly independent over $\Cl_n$.
To this end, suppose that
\ben
\sum_{\be,\Tc,\Uc}\ts c^{\be}_{\Tc\Uc}\tss\ze^{\be}_{\Tc\Uc}=0
\een
for some coefficients $c^{\be}_{\Tc\Uc}\in\Cl_n$. Let us act by the linear combination
on the left hand side on the basis vector $v^{}_{\ol\Vc}\in\wh U$ associated with
a standard barred tableau $\ol\Vc$ of a certain shape $\mu$. Let
$\ga$ denote the set of boxes of $\mu$ occupied by the barred entries
and let $\Vc$ be the tableau obtained from $\ol\Vc$ by removing all bars.
Using \eqref{idema} and the second expression
for $\ze^{\be}_{\Tc\Uc}$ in \eqref{zeeqet}, we get
\ben
\sum_{\be,\Tc}\ts c^{\be}_{\Tc\Vc}\tss
\phi^{}_{d(\Tc)}\ts e^{}_{\Rc^{\la,\be}}\ts \phi^*_{d(\Vc)}\ts v^{}_{\ol\Vc}=0.
\een
By applying \eqref{phiaa} we derive that
$
\phi^*_{d(\Vc)}\ts v^{}_{\ol\Vc}=b^*_{\Vc}\ts v^{}_{\Rc^{\mu,\ga}},
$
where $b^*_{\Vc}$ is a nonzero constant; cf. \cite[Lemma~7.8]{ww:ls}.
Therefore, by using \eqref{idema} again we conclude that
\ben
\sum_{\Tc}\ts c^{\ga}_{\Tc\Vc}\tss b^*_{\Vc}\tss
\phi^{}_{d(\Tc)} v^{}_{\Rc^{\mu,\ga}}=0,
\een
summed over standard $\mu$-tableaux $\Tc$. However, \eqref{phiaa} implies that
$\phi^{}_{d(\Tc)} v^{}_{\Rc^{\mu,\ga}}=b^{}_{\Tc}\tss v^{}_{\ol\Tc}$
for a nonzero constant $b^{}_{\Tc}$,
yielding
\ben
\sum_{\Tc}\ts c^{\ga}_{\Tc\Vc}\tss b^*_{\Vc}\tss b^{}_{\Tc} v^{}_{\ol\Tc}=0.
\een
Hence, using the definition \eqref{vbar} of $v^{}_{\ol\Tc}$ we get
$c^{\ga}_{\Tc\Vc}=0$ because the vectors $v^{}_{\Tc}$ are linearly
independent over $\Cl_n$. Thus, all elements $\ze^{\be}_{\Tc\Uc}$
are linearly independent over $\Cl_n$, as required.
\epf

\bre\label{rem:murphy}
The reader will have noticed that our logic in the arguments leading to the construction
of the basis of $\Sc_n$ in Theorem~\ref{thm:basalg} is different from
that of Murphy's paper \cite{m:rt}
devoted to the Hecke algebra.
We have taken advantage of the direct construction of the basis of the
$\Sc_n$-module $\wh U^{\la}$ given in \cite[Sec.~5.1]{hks:da} and \cite[Sec.~7.5]{ww:ls}.
On the other hand, it should be possible to reverse
the direction of the arguments for the Sergeev superalgebra
to get an analogue of the Murphy basis in the sense of \cite[Thm~3.9]{m:rt}; see also
\cite[Ch.~3]{m:ih}.
\qed
\ere

We will now combine Theorems~\ref{thm:primed} and \ref{thm:basalg} to give an
explicit realization of the irreducible $\Sc_n$-module $U^{\la}$ for any given $\la\Vdash n$.
We will produce
basis vectors of the left ideal $\Sc_n\tss e^{(r)}_{\Uc}$ and describe
the action of the generators of $\Sc_n$.
Since all left ideals $\Sc_n\tss e^{(r)}_{\Uc}$ associated with standard
barred tableaux of a given shape $\la$ are isomorphic as $\Sc_n$-modules, we choose
$\Uc$ to be the row-tableau $\Rc^{\la}$ with all entries unbarred. Furthermore,
we let $\Ec$ be the idempotent
in both Clifford superalgebras \eqref{cla} given by
\beql{ideme}
\Ec=\prod_{a=1}^m\frac{1+i\tss c_{2a-1}c_{2a}}{2}.
\eeq
If $\Tc$ is a standard $\la$-tableau with the diagonal entries $d_1<\dots<d_{\ell}$,
we let $\Ec^{\Tc}$ denote the corresponding idempotent defined as in \eqref{ideme}
for the Clifford superalgebra
with generators $c_{d_1},\dots,c_{d_{\ell}}$, where $m=\lfloor \ell/2\rfloor$.
We set $\Ec^{\la}=\Ec^{\Rc^{\la}}$.

By Theorem~\ref{thm:basalg}, the left ideal $J^{\la}:=\Sc_n\tss e^{}_{\Rc^{\la}}\Ec^{\la}$
is spanned over $\Cl_n$ by the elements
\beql{basele}
\ze^{\be}_{\Tc\Uc}\tss e^{}_{\Rc^{\la}}\Ec^{\la}.
\eeq
However, the second expression in
\eqref{zeequi} implies that
element \eqref{basele} is zero unless $\Uc=\Rc^{\la}$ and $\be$ is empty.
Due to Lemma~\ref{lem:inide} and property \eqref{phicli},
in this case the element can be written as
\beql{zedef}
\ze^{}_{\Tc}:=\phi^{}_{d(\Tc)}\ts e^{}_{\Rc^{\la}}\ts\Ec^{\la}=
e^{}_{\Tc}\ts\phi^{}_{d(\Tc)}\ts \Ec^{\la}=
e^{}_{\Tc}\ts\Ec^{\Tc}\phi^{}_{d(\Tc)}=\Ec^{\Tc} e^{}_{\Tc}\ts\phi^{}_{d(\Tc)}.
\eeq
Hence, $J^{\la}$ is spanned over $\Cl_n$ by the vectors $\ze^{}_{\Tc}$ associated with
the standard $\la$-tableaux $\Tc$.

Furthermore, for a fixed standard $\la$-tableau $\Tc$, the span $\Cl_n\tss \Ec^{\Tc}$
is a free module over the Clifford superalgebra $\Cl^{\Tc}_{n-m}\subset\Cl_n$.
The latter is generated by the elements $c_{d_1},c_{d_3},\dots,c_{d_{2m-1}}$ together with
the $c_a$ for $a$ running over the set $\{1,\dots,n\}\setminus\{d_1,d_2,\dots,d_{2m}\}$.
The left ideal $J^{\la}$ thus has a basis over $\CC$ arising from the decomposition
\beql{lefde}
J^{\la}=\bigoplus_{\sh(\Tc)=\la}\ts \Cl^{\Tc}_{n-m}\tss \ze^{}_{\Tc}.
\eeq
This is clear from the observation that
the number of spanning vectors coincides with the dimension $2^{n-m}g_{\la}$
of the irreducible module $U^{\la}$.
We will now derive explicit formulas for the action of the generators of $\Sc_n$
in this basis. They will essentially coincide with \eqref{saavt} and \eqref{jmacti}, but the latter
were applied to a basis of the module
$\wh U^{\la}$ isomorphic to a direct sum of several copies of $U^{\la}$;
the identification of $U^{\la}$ inside $\wh U^{\la}$ required some additional steps;
cf. \cite[Thm~5.1.3]{hks:da}, \cite[Thm~8.3]{n:ys} and \cite[Thm~7.17]{ww:ls}.
It will be convenient to normalize the basis vectors by setting
\ben
\xi^{}_{\Tc}=\frac{1}{b^{}_{\Tc}}\ts\ze^{}_{\Tc}
\een
for every standard $\la$-tableau $\Tc$, where the nonzero constant $b^{}_{\Tc}$
is defined as in the proof of Theorem~\ref{thm:basalg} by the relation
$\phi^{}_{d(\Tc)} v^{}_{\Rc^{\la}}=b^{}_{\Tc}\tss v^{}_{\Tc}$.
We use notation \eqref{yatc}.

\bco\label{cor:baact}
The generators of the Sergeev superalgebra $\Sc_n$ act in the normalized basis of the simple
module $J^{\la}$ defined by the decomposition \eqref{lefde} by the rule:
\begin{align}
s_a\tss \xi^{}_{\Tc}&=\Bigg(\frac{1}{\ka_{a+1}(\Tc)-\ka_{a}(\Tc)}
+\frac{c_ac_{a+1}}{\ka_{a+1}(\Tc)+\ka_{a}(\Tc)}\Bigg)\tss \xi^{}_{\Tc}
+\Yc_a(\Tc)\ts\xi^{}_{s_a\Tc},
\label{saxit}\\[0.5em]
x_a\tss \xi^{}_{\Tc}&=\ka_{a}(\Tc)\ts \xi^{}_{\Tc},
\label{xaxi}
\end{align}
where we assume that $\xi^{}_{s_a\Tc}=0$ if the tableau $s_a\Tc$ is not standard. Moreover,
for a given $\Tc$, the action of the elements $c_{d_2},c_{d_4},\dots,c_{d_{2m}}\in\Cl_n$ is determined by
\beql{clac}
c_{d_{2a}} \xi^{}_{\Tc}=i\ts c_{d_{2a-1}} \xi^{}_{\Tc},\qquad a=1,\dots,m.
\eeq
\eco

\bpf
Relation~\eqref{xaxi} follows from
\eqref{xietsign} and the definition \eqref{zedef} of $\ze^{}_{\Tc}$.
To verify \eqref{saxit}, suppose first that the tableau $s_a\Tc$ is standard and
the length of the permutation $d(s_a\Tc)$ is greater than
the length of $d(\Tc)$.
Then the definition of $\ze^{}_{\Tc}$
yields $\phi_a\tss \ze^{}_{\Tc}=\ze^{}_{s_a\Tc}$. Using \eqref{phiadef} and \eqref{xaxi}
we derive that
\ben
s_a\tss \ze^{}_{\Tc}=\Bigg(\frac{1}{\ka_{a+1}(\Tc)-\ka_{a}(\Tc)}
+\frac{c_ac_{a+1}}{\ka_{a+1}(\Tc)+\ka_{a}(\Tc)}\Bigg)\tss \ze^{}_{\Tc}
+\frac{1}{\ka_{a}(\Tc)^2-\ka_{a+1}(\Tc)^2}\ts\ze^{}_{s_a\Tc}.
\een
Relation \eqref{saxit} now follows from
the identity for the normalizing constants
\ben
b_{s_a\Tc}=b_{\Tc}\ts \big(\ka_{a}(\Tc)^2-\ka_{a+1}(\Tc)^2\big)\ts \Yc_a(\Tc)
\een
implied by \eqref{phiaa}. Furthermore,
we also have the relation $\phi^2_a\tss \ze^{}_{\Tc}=\phi_a\tss\ze^{}_{s_a\Tc}$,
whose left hand side can be evaluated with the use of \eqref{phiasqa} and \eqref{xaxi}.
Calculating as above by swapping the roles of $\Tc$ and $s_a\Tc$, we get \eqref{saxit}
in the case, where the length of the permutation $d(s_a\Tc)$
is less than that
of $d(\Tc)$. If the tableau $s_a\Tc$ is not
standard, then $\phi_a\tss \ze^{}_{\Tc}=0$ by Lemma~\ref{lem:inide},
thus implying \eqref{saxit} in this case.
Finally, relation \eqref{clac} is immediate from \eqref{ideme}, since
$(c_{2a}-i\ts c_{2a-1})\ts\Ec=0$ for $a=1,\dots,m$.
\epf

The above arguments show that the isomorphism
\eqref{isomsn} with $\Uc=\Rc^{\la}$ can be described explicitly by
\beql{vtmap}
v^{}_{\Tc}\mapsto \frac{1}{b_{\Tc}}\ts\phi^{}_{d(\Tc)}\tss e^{}_{\Rc^{\la}}
=\frac{1}{b_{\Tc}}\ts e^{}_{\Tc}\tss\phi^{}_{d(\Tc)}.
\eeq
Moreover, by Corollary~\ref{cor:baact}, the simple $\Sc_n$-module $U^{\la}$
is isomorphic to the quotient of $\wh U^{\la}$ by the submodule generated by
all elements of the form $(c_{d_{2a}}-i\ts c_{d_{2a-1}})\tss v^{}_{\Tc}$
with $a=1,\dots,m$ and $\Tc$ running over the standard $\la$-tableaux.

\bre\label{rem:affine}
The formulas of Corollary~\ref{cor:baact} also provide a realization of
the irreducible {\em calibrated} (or {\em completely splittable})
module over the {\em degenerate affine Hecke--Clifford algebra} $\Sc^{\text{\rm aff}}_n$
associated with the strict partition $\la$; see \cite{hks:da} and \cite{w:cs}
for the classification of the irreducible modules of this class.
The elements $x_a$ should be understood as generators of $\Sc^{\text{\rm aff}}_n$; see \cite{n:ys}.
\ere

\subsection{Seminormal form for the spin symmetric group algebra}
\label{subsec:spin}

Consider the restriction of the representation $\wh U^{\la}$ of the Sergeev superalgebra
defined in \eqref{ula} to
the spin symmetric group algebra
$\CC\Sym_n^-$. Rewriting \eqref{embsy}, we get
\beql{embsyinv}
t_a=y_a\tss s_a,\qquad y_a:=\frac{c_{a+1}-c_a}{\sqrt2}.
\eeq
Hence, by deriving from \eqref{saavt}, we get
\beql{taspin}
t_a\tss v^{}_{\Tc}=
\frac{\sqrt2\ts\big(\ka_{a}(\Tc)\tss c^{}_{a}-\ka_{a+1}(\Tc)\tss c^{}_{a+1}\big)}
{\ka_{a}(\Tc)^2-\ka_{a+1}(\Tc)^2}
\tss v^{}_{\Tc}
+y_a\ts \Yc_a(\Tc)\ts v^{}_{s_a\Tc}.
\eeq
We will use these formulas to calculate the action of the generators $t_a$
on the vectors $\theta^{}_{\Tc}\in \wh U^{\la}$ defined by
\beql{thetadef}
\theta^{}_{\Tc}=y_{a_r}\dots y_{a_1}\ts v^{}_{\Tc}
\eeq
with the use of a reduced decomposition $d(\Tc)=s_{a_1}\dots s_{a_r}$. The expression
in \eqref{thetadef} depends on the decomposition
up to a sign, due to the relations satisfied by
the elements $y_a$:
\beql{ybra}
y_a^2=-1,\qquad y_ay_{a+1}y_a=y_{a+1}y_a y_{a+1},\qquad y_ay_b=-y_by_a,\qquad |a-b|>1.
\eeq
For this reason, we will fix particular reduced decompositions of the permutations
$d(\Tc)$ as given below. Any different choice would only affect the signs $\ve_a(\Tc)$
appearing in Lemma~\ref{lem:spin}.

Consider the {\em row-word} $\row(\Tc)$ of any $\la$-tableau $\Tc$
as the sequence of the entries of $\Tc$ read consecutively
by rows from left to right, starting from the top row.
For every $b=2,\dots,n$, the $b$-{\em inversion} in $\row(\Tc)$ is a pair $(b',b)$
such that $b'<b$ and $b'$ occurs after $b$ in $\row(\Tc)$. Denote by $i_b=i_b(\Tc)$ the number
of $b$-inversions in $\row(\Tc)$. Then the formula
\beql{dtre}
d(\Tc)=\prod_{b=2,\dots,n}^{\longrightarrow}\ts (s^{}_{b-1}s^{}_{b-2}\dots s^{}_{b-i_b})
\eeq
defines a reduced decomposition; if $i_b=0$, then the corresponding
product term is considered to be equal to $1$. We will use \eqref{dtre} to define
the vectors $\theta^{}_{\Tc}$ in \eqref{thetadef}. For all $b=1,\dots,n$ set
$b^{\Tc}:=d(\Tc)^{-1}(b)$.

\ble\label{lem:spin}
The generators of the superalgebra $\CC\Sym_n^-$ act on the vectors $\theta^{}_{\Tc}\in \wh U^{\la}$
by the rule:
\beql{saxitspin}
t_a\tss \theta^{}_{\Tc}=
\frac{\sqrt2\ts\big(\ka_{a}(\Tc)\tss c^{}_{a^{\Tc}}-\ka_{a+1}(\Tc)\tss c^{}_{(a+1)^{\Tc}}\big)}
{\ka_{a}(\Tc)^2-\ka_{a+1}(\Tc)^2}
\tss \theta^{}_{\Tc}
+\ve_a(\Tc)\tss \Yc_a(\Tc)\ts\theta^{}_{s_a\Tc},
\eeq
where
we assume that $\theta^{}_{s_a\Tc}=0$ if the tableau $s_a\Tc$ is not standard, while
\ben
\ve_a(\Tc)=(-1)^{i_a(i_{a+1}+1)+i_{a+2}+\dots+i_n}\qquad\text{or}\qquad
\ve_a(\Tc)=(-1)^{(i_a+1)(i_{a+1}+1)+i_{a+2}+\dots+i_n},
\een
depending on whether the length of $d(s_a\Tc)$ is greater or smaller than
that of $d(\Tc)$.
\ele

\bpf
Using \eqref{thetadef} write
\ben
t_a\tss \theta^{}_{\Tc}=(-1)^r y_{a_r}\dots y_{a_1}\ts t_a\tss  v^{}_{\Tc}
\een
and apply \eqref{taspin} to expand $t_a\tss  v^{}_{\Tc}$.  Since $y_a^2=-1$,
by \eqref{thetadef} we have
\ben
v^{}_{\Tc}=(-1)^r y_{a_1}\dots y_{a_r}\ts\theta^{}_{\Tc}.
\een
Observe that
\ben
y_a c_ay_a=c_{a+1},\qquad y_a c_{a+1}y_a=c_a,\qquad y_a c_b\tss y_a=c_b,\quad b\ne a,a+1.
\een
Hence, for any $b\in\{1,\dots,n\}$ we have
\ben
y_{a_r}\dots y_{a_1}c_b\ts y_{a_1}\dots y_{a_r}=c^{}_{b^{\Tc}}
\een
which is independent of the chosen reduced decomposition of $d(\Tc)$. This verifies the
coefficient of $\theta^{}_{\Tc}$ on the right hand side of \eqref{saxitspin}.

Now suppose that the tableau $s_a\Tc$ is standard and
the length of $d(s_a\Tc)$ is greater than
that of $d(\Tc)$. If
\ben
d(s_a\Tc)=s_{b_1}\dots s_{b_{r+1}}
\een
is the reduced decomposition chosen by the above rule, then similar to
the previous argument, the factor
$\ve_a(\Tc)$ appearing in \eqref{saxitspin} will be given by
\beql{epst}
\ve_a(\Tc)=-y_{a_r}\dots y_{a_1}y_a\tss y^{}_{b_1}\dots y^{}_{b_{r+1}}.
\eeq
In this case, the row-words have the form
\ben
\row(\Tc)=(1,\dots,a,\dots,a+1,\dots)\Fand \row(s_a\Tc)=(1,\dots,a+1,\dots,a,\dots)
\een
so that $\row(s_a\Tc)$ is obtained from $\row(\Tc)$ by swapping $a$ and $a+1$. Therefore,
we have
$i_a(s_a\Tc)=i_{a+1}$ and $i_{a+1}(s_a\Tc)=i_{a}+1$, whereas $i_b(s_a\Tc)=i_{b}$
for $b\ne a,a+1$. Hence, for the part of the product in \eqref{epst} we can write
\begin{multline}
\non
y_a\tss y^{}_{b_1}\dots y^{}_{b_{r+1}}=
y_a\prod_{b=2,\dots,a-1}^{\longrightarrow}\ts (y^{}_{b-1}\dots y^{}_{b-i_b})\\[-0.3em]
{}\times\tss(y^{}_{a-1}\dots y^{}_{a-i_{a+1}})(y^{}_{a}\dots y^{}_{a-i_{a}})
\prod_{b=a+2,\dots,n}^{\longrightarrow}\ts (y^{}_{b-1}\dots y^{}_{b-i_b}).
\end{multline}
By using properties \eqref{ybra}, we will move $y_a$ to the right through the first product, taking the
signs into account, then apply the identity
\ben
y_a\tss (y^{}_{a-1}\dots y^{}_{a-i_{a+1}})(y^{}_{a}\dots y^{}_{a-i_{a}})
=(-1)^{(i_a+1)i_{a+1}+1}(y^{}_{a-1}\dots y^{}_{a-i_{a}})(y^{}_{a}\dots y^{}_{a+1-i_{a+1}}).
\een
To verify the latter, first note the relation
\beql{ykmo}
y_k\tss (y^{}_{a}\dots y^{}_{a-i_{a}})=(-1)^{i_a+1}(y^{}_{a}\dots y^{}_{a-i_{a}})\tss y_{k+1}
\eeq
which holds for any $k\in\{a-i_{a},\dots,a-1\}$. Indeed, to check \eqref{ykmo}, move $y_k$ to the right
by permuting it with $y^{}_{a},\dots,y^{}_{k+2}$, then apply
the braid relation $y^{}_ky^{}_{k+1}y^{}_k=y^{}_{k+1}y^{}_ky^{}_{k+1}$,
and finally move $y_{k+1}$ to the right
by permuting it with $y^{}_{k-1},\dots,y^{}_{a-i_{a}}$. The identity follows
by applying \eqref{ykmo} consecutively for $k=a-i_{a+1},\dots,a-1$ and taking into account
that $y_a^2=-1$.

The calculation thus yields the formula
\ben
y_a\tss y^{}_{b_1}\dots y^{}_{b_{r+1}}=(-1)^{i_2+\dots+i_{a-1}+(i_a+1)i_{a+1}+1}\ts y_{a_1}\dots y_{a_r}.
\een
Hence, \eqref{epst} implies the required expression
\ben
\ve_a(\Tc)=-y_{a_r}\dots y_{a_1}y_a\tss y^{}_{b_1}\dots y^{}_{b_{r+1}}
=(-1)^{i_2+\dots+i_{a-1}+(i_a+1)i_{a+1}+r}
=(-1)^{i_a(i_{a+1}+1)+i_{a+2}+\dots+i_n}
\een
in the case under consideration, since $r=i_2+\dots+i_n$.
By using the relation $t_a^2=1$ and applying $t_a$ to
both sides of \eqref{saxitspin} we derive that
$\ve_a(\Tc)\ve_a(s_a\Tc)=1$. This gives the required formula in
the remaining case.
\epf

Observe that the right hand side of
\eqref{saxitspin} only depends  on the Clifford generators,
corresponding to the non-diagonal entries of the row-tableau $\Rc^{\la}$.
Indeed, if $b^{\Tc}=d(\Tc)^{-1}(b)$ is a diagonal entry of $\Rc^{\la}$, then
$b=d(\Tc)(b^{\Tc})$ is a diagonal entry of $\Tc$. However, in this case,
$\ka_{b}(\Tc)=0$.

Denote by $\Cl^{\la}_{n-\ell}$ the subalgebra of the Clifford superalgebra
$\Cl_n$ generated by the elements
$c_{b_1},\dots,c_{b_{n-\ell}}$,
where $b_1<\dots<b_{n-\ell}$ are the non-diagonal entries of $\Rc^{\la}$.
The superalgebra $\Cl^{\la}_{n-\ell}$
has a unique right module $E^{\la}$ of dimension $2^p$ with $p=\lceil \frac{n-\ell}{2}\rceil$.
We will realize the basis vectors of $E^{\la}$ inside the Clifford superalgebra
$\Cl_n$ as follows.
First consider the Clifford superalgebra $\Cl_{2p}$ with its elements
\ben
\Fc^{\ts+}_a=1+i\tss c_{2a-1}c_{2a}\Fand \Fc^{\ts-}_a =c_{2a-1}+i\tss c_{2a},
\qquad a=1,\dots,p,
\een
and note the relations
\ben
\Fc^{\ts+}_a\tss c_{2a-1}=\Fc^{\ts-}_a,\qquad
\Fc^{\ts+}_a\tss c_{2a}=-i\ts \Fc^{\ts-}_a,\qquad
\Fc^{\ts-}_a\tss c_{2a-1} =-\Fc^{\ts+}_a,\qquad \Fc^{\ts-}_a\tss c_{2a} =-i\ts \Fc^{\ts+}_a.
\een
The right $\Cl_{2p}$-module $E$ has the basis consisting of the $2^p$ elements
\ben
\Fc^{\ts\de}=\Fc_1^{\ts\de_1}\dots \Fc_p^{\ts\de_p},
\een
where $\de=(\de_1,\dots,\de_p)$ with each $\de_a$ taking the value $+$ or $-$. Similarly,
the right $\Cl_{2p+1}$-module $E$ has the basis consisting of the $2^{p+1}$ elements
\ben
\Fc^{\ts\de}=\Fc_1^{\ts\de_1}\dots \Fc_p^{\ts\de_p}\ts c_{2p+1}^{\al},
\een
where we now set $\de=(\de_1,\dots,\de_p,\al)$
with the same parameters $\de_a$ as above and $\al$ taking two possible values $0$ and $1$.

Returning to the Clifford superalgebra $\Cl^{\la}_{n-\ell}$ with $n-\ell=2p$ or $n-\ell=2p+1$,
we let $E^{\la}$ denote the corresponding right $\Cl^{\la}_{n-\ell}$-module constructed
as the module $E$ over $\Cl_{2p}$ or $\Cl_{2p+1}$ in the same way above.
We will keep the notation
$\Fc^{\ts\de}$ for the corresponding basis vectors of $E^{\la}$.

\bth\label{thm:spin}
The vectors $\Fc^{\ts\de}\tss \theta^{}_{\Tc}$ parameterized by the tuples $\de$
and standard $\la$-tableaux $\Tc$ form a basis of the irreducible representation of $\CC\Sym_n^-$
associated with $\la$.
The action of the generators $t_a$ of $\CC\Sym_n^-$ in the basis is determined by
the formulas of Lemma~\ref{lem:spin}.
\eth

\bpf
The span of the vectors is invariant under the action of $\CC\Sym_n^-$. Their number
equals $2^r g_{\la}$ with $r=\lceil \frac{n-\ell}{2}\rceil$ which coincides with the dimension
of the simple $\CC\Sym_n^-$-module $V^{\la}$ associated with the shifted diagram $\la$.
\epf

\bex\label{ex:threeone}
To illustrate Theorem~\ref{thm:spin}, consider the case $n=4$ with $\la=(3,1)$.
The row-tableau is found by
\ben
\Rc^{(3,1)}\ \ =\quad\begin{ytableau}
    1 & 2 & 3\\
    \none&4
\end{ytableau}
\een
so that the Clifford superalgebra $\Cl^{\la}_{n-\ell}=\Cl^{(3,1)}_2$
is generated by the elements $c_2$ and $c_3$.
The simple module $V^{(3,1)}$ is four-dimensional and we have
the basis vectors
\ben
(1+i\tss c_{2}c_{3})\ts \theta_{\Tc},\qquad (1+i\tss c_{2}c_{3})\ts \theta_{\Uc},\qquad
(c_{2}+i\tss c_{3})\ts \theta_{\Tc},\qquad (c_{2}+i\tss c_{3})\ts \theta_{\Uc},
\een
for the standard tableaux
\ben
\Tc\ \ =\quad\begin{ytableau}
    1 & 2 & 3\\
    \none&4
\end{ytableau}
\Fand
\Uc\ \ =\quad\begin{ytableau}
    1 & 2 & 4\\
    \none&3
\end{ytableau}\ \ .
\een
By Lemma~\ref{lem:spin}, we have
\ben
t_1\ts \theta_{\Tc}=c_2\ts\theta_{\Tc},\qquad
t_2\ts \theta_{\Tc}=\frac{-c_2+\sqrt3\ts c_3}{2}\ts\theta_{\Tc},\qquad
t_3\ts \theta_{\Tc}=\frac{\sqrt3}{3}\ts c_3\ts \theta_{\Tc}+\frac{\sqrt6}{3}\ts \theta_{\Uc}
\een
and
\ben
t_1\ts \theta_{\Uc}=c_2\ts\theta_{\Uc},\qquad
t_2\ts \theta_{\Uc}=c_2\ts\theta_{\Uc},\qquad
t_3\ts \theta_{\Uc}=\frac{\sqrt3}{3}\ts c_3\ts \theta_{\Uc}+\frac{\sqrt6}{3}\ts \theta_{\Tc}.
\een
This determines the action of the generators on all basis vectors of $V^{(3,1)}$.
For instance,
\ben
t_3\ts (1+i\tss c_{2}c_{3})\ts \theta_{\Tc}=
(1+i\tss c_{2}c_{3})\ts t_3\ts \theta_{\Tc}=
(1+i\tss c_{2}c_{3})\ts
\Big(\frac{\sqrt3}{3}\ts c_3\ts \theta_{\Tc}+\frac{\sqrt6}{3}\ts \theta_{\Uc}\Big)
\een
which equals
\ben
-\frac{i\tss\sqrt3}{3}\ts (c_{2}+i\tss c_{3})\ts \theta_{\Tc}
+\frac{\sqrt6}{3}\ts (1+i\tss c_{2}c_{3})\ts \theta_{\Uc}.
\een
\eex

\section{Fusion procedure}
\label{sec:fp}

Following \cite{n:ys, n:ci},
introduce rational functions in variables $u,v$ with values in $\Sc_n$
by setting
\ben
\vp_{ab}(u,v)=1+t_{ab}\ts \frac{\sqrt{2}\ts(u\tss c_a-v\tss c_b)}{u^2-v^2},\qquad a\ne b.
\een
Equivalently,
in terms of the presentation $\Sc_n=\CC\Sym_n\ltimes\Cl_n$ as defined in \eqref{symgr},
we have
\ben
\vp_{ab}(u,v)=1-\frac{(a,b)}{u-v}+\frac{(a,b)\ts c_ac_b}{u+v}.
\een
Take $n$ complex variables $u_1,\dots,u_n$
and set
\beql{phiuuu}
\Phi(u_1,\dots,u_n)=
\prod_{1\leqslant a<b\leqslant n}\vp^{}_{ab}(u_a,u_b),
\eeq
where the product is taken in the lexicographic order on the pairs $(a,b)$.

Suppose that $\la\Vdash n$
and let $\Uc$ be a standard barred $\la$-tableau.
For every $a=1,\dots,n$ set
$\ka_a=\ka_{\asf}(\Uc)$ if $\asf=a$ or $\asf=\bar a$ is the entry of $\Uc$.

\bth\label{thm:fusion}
The consecutive evaluations are well-defined
and we have the identity
\beql{evide}
\Phi(u_1,\dots,u_n)\big|_{u_1=\ka_{\mathsf 1}}\big|_{u_2=\ka_2}\dots \big|_{u_n=\ka_n}
=\frac{n!}{g_{\la}}\ts e^{}_{\Uc}.
\eeq
\eth

\bpf
We argue by induction on $n$ and note that the identity is trivial for $n=1$.
Observe that $\vp_{ab}(u,v)$ and $\vp_{cd}(u',v')$ commute if the indices $a,b,c,d$ are
distinct. Hence
using
the induction hypothesis and setting $u=u_n$, for $n\geqslant 2$ we get
\ben
\Phi(u_1,\dots,u_n)\big|_{u_1=\ka_1}\dots \big|_{u_{n-1}=\ka_{n-1}}
=\frac{(n-1)!}{g_{\mu}}\ts e^{}_{\Vc}\ts
\vp^{}_{1n}(\ka_1,u)\dots \vp^{}_{n-1,n}(\ka_{n-1},u),
\een
where the barred tableau $\Vc$ is obtained from $\Uc$ by
removing the box occupied by $n$ or $\bar n$ and $\mu$ is the shape of $\Vc$.
Hence, it suffices to show that the following evaluation is well-defined with the value
given by
\beql{eva}
e^{}_{\Vc}\ts
\vp^{}_{1n}(\ka_1,u)\dots \vp^{}_{n-1,n}(\ka_{n-1},u)\big|_{u=\ka_n}
=\frac{n\tss g_{\mu}}{g_{\la}}\ts e^{}_{\Uc}.
\eeq

\ble\label{lem:prjm}
We have the relation
\beql{rede}
e^{}_{\Vc}\ts \vp^{}_{n,n-1}(u,\ka_{n-1})\dots \vp^{}_{n1}(u,\ka_{1})
=\frac{u-x_n}{u}\ts e^{}_{\Vc}.
\eeq
\ele

\bpf
A counterpart of
this relation was pointed out in \cite[Prop.~3.1]{n:ci} in a different
context; we will give a direct proof.
It will be convenient to prove more general relations in $\Sc_n$,
\beql{redegen}
e^{}_{\Vc}\ts \vp^{}_{n,r}(u,\ka_{r})\dots \vp^{}_{n1}(u,\ka_{1})
=\frac{u-x_{r n}}{u}\ts e^{}_{\Vc},
\eeq
where $r=1,\dots,n-1$ and $\Vc$ now denotes a standard barred tableau
with $r$ boxes; the numbers $\ka_a=\ka_{\asf}(\Vc)$ are the signed contents of the entries of $\Vc$
and
\ben
x_{rn}=\sqrt{2}\ts (t_{1n}+\dots+t_{rn})\ts c_n.
\een
Note that $x_{rn}$ commutes with $x_1,\dots,x_r$ and so it commutes with $e^{}_{\Vc}$.
It is clear that \eqref{redegen} holds for $r=1$. To proceed by
the induction on $r$, suppose that $r\geqslant 2$.
Let $\Wc$ be the barred tableau obtained from $\Vc$ by
removing the box occupied by $r$ or $\bar r$. Write $e^{}_{\Vc}=e^{}_{\Vc}e^{}_{\Wc}$,
which is an obvious generalization of the idempotent property established in
Proposition~\ref{prop:idemp},
and observe that $e^{}_{\Wc}$ commutes with the factor $\vp^{}_{n,r}(u,\ka_{r})$. Hence,
by the induction hypothesis, the left hand side of \eqref{redegen}
can be written as
\ben
e^{}_{\Vc}\ts \vp^{}_{n,r}(u,\ka_{r})\ts \frac{u-x_{r-1,n}}{u}\ts e^{}_{\Wc}
=e^{}_{\Vc}\ts \Big(1+t_{rn}\ts \frac{\sqrt{2}\ts(\ka_rc_r-uc_n)}{u^2-\ka_r^2}\Big)
\big(1-\frac{x_{r-1,n}}{u}\big).
\een
To simplify
this expression, use \eqref{xietsign} to derive
the relation
\begin{multline}
e^{}_{\Vc}\ts t_{rn}\ts (t_{1n}+\dots+t_{r-1,n})
=-e^{}_{\Vc}\ts (t_{1r}+\dots+t_{r-1,r})\ts t_{rn}\\
=\frac{1}{\sqrt{2}}\ts
e^{}_{\Vc}\ts x_r\ts c_r\ts t_{rn}=-\frac{1}{\sqrt{2}}\ts
\ka_r\ts e^{}_{\Vc}\ts t_{rn}\ts c_r
\non
\end{multline}
thus verifying \eqref{redegen}.
\epf

Returning to the proof of the theorem,
use the relations
\ben
\vp_{ab}(u,v)\ts \vp_{ba}(v,u)=A(u,v)
\een
which hold for $a\ne b$,
where $A(u,v)$ is defined in \eqref{auv}.
Hence, we derive from Lemma~\ref{lem:prjm} that
\ben
e^{}_{\Vc}\ts
\vp^{}_{1n}(\ka_1,u)\dots \vp^{}_{n-1,n}(\ka_{n-1},u)
=\frac{u}{u-x_n}\ts e^{}_{\Vc}\ts A(u,\ka_1)\dots A(u,\ka_{n-1}),
\een
and the right hand side can be written as
\ben
\frac{u-\ka_n}{u-x_n}\ts e^{}_{\Vc}\cdot
\frac{u}{u-\ka_n}\ts A(u,\ka_1)\dots A(u,\ka_{n-1}).
\een
The required relation \eqref{eva} will follow from the claim that both the following
evaluations are well-defined with the values given by
\beql{ratfu}
\frac{u-\ka_n}{u-x_n}\ts e^{}_{\Vc}\big|_{u=\ka_n}=e^{}_{\Uc}
\eeq
and
\beql{ratfnu}
\frac{u}{u-\ka_n}\ts A(u,\ka_1)\dots A(u,\ka_{n-1})\big|_{u=\ka_n}=\frac{n\tss g_{\mu}}{g_{\la}}.
\eeq
To verify \eqref{ratfu}, note that the definition \eqref{murphyse} implies
\beql{eutet}
e^{}_{\Vc}=\sum_{\Vc\nearrow \Wc} e^{}_{\Wc},
\eeq
where $\Vc\sclnearrow \Wc$ means that
the barred tableau $\Wc$ is obtained from $\Vc$
by adding one box with the entry $n$ or $\bar n$.
Since $x_n e^{}_{\Wc}=\ka_{\nsf}(\Wc)e^{}_{\Wc}$ for $\nsf=n$ or $\nsf=\bar n$
by \eqref{xietsign}, we use \eqref{xnid}
to write
\ben
e^{}_{\Vc}\ts \frac{u-\ka_n\ }{u-x_n}=\sum_{\Vc\nearrow \Wc}
e^{}_{\Wc}\ts\frac{u-\ka_n\ }{u-\ka_{\nsf}(\Wc)}
=e^{}_{\Uc}+\sum_{\Vc\nearrow \Wc,\ts \Wc\ne\Uc}
e^{}_{\Wc}\ts\frac{u-\ka_n\ }{u-\ka_{\nsf}(\Wc)}.
\een
Hence the evaluation at $u=\ka_n$ is well-defined and \eqref{ratfu} follows.

To verify \eqref{ratfnu}, note that the product
depends only on the shape $\mu$ of $\Vc$ so we may choose a particular
tableau $\Vc$ for its evaluation. Moreover, since the rational function $A(u,v)$
is even in $u$ and $v$, we may
take $\Vc$ to be unbarred; we will take
the row tableau
obtained by filling in the boxes of $\mu$ with
the numbers $1,\dots,n-1$ by consecutive rows from left to right
in each row.
Noting the formula
\ben
A(u,0)A(u,\sqrt 2)\dots A\big(u,\sqrt{p(p-1)}\ts\big)=\frac{u^2-p(p+1)}{u^2-p(p-1)},
\een
we can write
\ben
A(u,\ka_1)\dots A(u,\ka_{n-1})=\prod_{a=1}^m\ts \frac{u^2-\mu_a(\mu_a+1)}{u^2-\mu_a(\mu_a-1)},
\een
where we set $\mu=(\mu_1,\dots,\mu_m)$ with $m=\ell(\mu)$. Suppose first, that
$m=\ell(\la)$ so that the tableau $\Uc$ is obtained from $\Vc$ by adding a
box with the entry $n$ or $\bar n$ at the row $k\in\{1,\dots,m\}$. Then
$\ka_n=\sqrt{\mu_k(\mu_k+1)}$ or $\ka_n=-\sqrt{\mu_k(\mu_k+1)}$, respectively,
while $\mu_{k-1}\geqslant \mu_k+2$. Therefore, the left hand side of \eqref{ratfnu}
is well-defined with the value
\ben
(\mu_k+1)\ts\prod_{a=1,\ a\ne k}^m\
\frac{\mu_k(\mu_k+1)-\mu_a(\mu_a+1)}{\mu_k(\mu_k+1)-\mu_a(\mu_a-1)}.
\een
It is easily seen that this
coincides with $n\ts g_{\mu}/g_{\la}$, as implied by the Schur formula \eqref{gla},
thus confirming \eqref{ratfnu} in the case under consideration.
In the remaining case with $m=\ell(\la)-1$, the tableau $\Uc$
is obtained from $\Vc$ by adding the
box $(\ell(\la),\ell(\la))$ with the entry $n$. This is possible when $\mu_m\geqslant 2$.
In this case, $\ka_n=0$ and
the left hand side of \eqref{ratfnu}
is well-defined with the value
\ben
\prod_{a=1}^m\ts \frac{\mu_a+1}{\mu_a-1},
\een
which coincides with the right hand side of \eqref{ratfnu}.
\epf

\bex\label{ex:twoone}
Consider the tableau
\ben
\Uc\ \ =\quad\begin{ytableau}
    1 & 2\\
    \none&3
\end{ytableau}
\een
of shape $(2,1)$. We have
\begin{multline}
\Phi(u_1,u_2,u_3)\\
{}=\big(1+t_{12}\ts \frac{(u_1c_1-u_2c_2)\sqrt{2}}{u_1^2-u_2^2}\big)
\big(1+t_{13}\ts \frac{(u_1c_1-u_3c_3)\sqrt{2}}{u_1^2-u_3^2}\big)
\big(1+t_{23}\ts \frac{(u_2c_2-u_3c_3)\sqrt{2}}{u_2^2-u_3^2}\big).
\non
\end{multline}
Note that the direct evaluation of this function at
$u_1=u_3=0$ and $u_2=\sqrt{2}$ is not defined. As in the proof of Theorem~\ref{thm:fusion},
we first set
$u_1=0$, $u_2=\sqrt{2}$ and write $u$ for $u_3$. Then
\ben
\Phi(0,\sqrt{2},u)
=\big(1+t_{12}c_2\big)
\big(1+\ts \frac{t_{13}\tss c_3\sqrt{2}}{u}\big)
\big(1+t_{23}\ts \frac{2c_2-uc_3\sqrt{2}}{2-u^2}\big).
\een
Now observe that $1+t_{12}c_2=(1+t_{12}c_2)\ts t_{12}c_2$ and
\ben
t_{12}c_2 \cdot \frac{t_{13}\tss c_3\sqrt{2}}{u}\cdot t_{23}\ts \frac{2c_2-uc_3\sqrt{2}}{2-u^2}
=\frac{t_{12}t_{13}t_{23}\ts (2uc_2+2c_3\sqrt{2})}{u(u^2-2)}
=\frac{t_{13}\ts (2uc_2+2c_3\sqrt{2})}{u(u^2-2)}.
\een
Hence
\ben
\Phi(0,\sqrt{2},u)
=\big(1+t_{12}c_2\big)\big(1-t_{13}\ts \frac{2c_2+uc_3\sqrt{2}}{2-u^2}
+t_{23}\ts \frac{2c_2-uc_3\sqrt{2}}{2-u^2}\big),
\een
which is regular at $u=0$ with the value
\ben
\big(1+t_{12}c_2\big)\big(1-(t_{13}-t_{23})c_2\big)
=\big(1+t_{12}c_2\big)\big(1+c_2(t_{13}-t_{23})\big).
\een
On the other hand, the corresponding idempotent from Example~\ref{ex:smdi3}
takes the form
\ben
e^{}_{\Uc}=\frac{\sqrt{2}+x_2}{2\tss\sqrt{2}}\cdot \frac{6-x_3^2}{6}.
\een
Since $x_2=\sqrt{2}\ts t_{12}c_2$ and
\ben
x_3^2=2\ts m_3^2=2(t_{13}+t_{23})^2=4+2\ts t_{12}(t_{13}-t_{23}),
\een
we get
\ben
e^{}_{\Uc}=\frac{1}{6}\ts\big(1+t_{12}c_2\big)\big(1-t_{12}(t_{13}-t_{23})\big)
=\frac{1}{6}\ts\big(1+t_{12}c_2\big)\big(1+c_2(t_{13}-t_{23})\big)
\een
thus agreeing with Theorem~\ref{thm:fusion}.
\eex

\bex\label{ex:symm}
Consider the particular case $\la=(n)$ with the only standard unbarred tableau $\Uc$.
This time the evaluations in \eqref{evide} can be made simultaneously, so that
the evaluated product is given by
\ben
\Phi\big(0,\sqrt{2},\dots,\sqrt{n(n-1)}\ts\big)=
\prod_{1\leqslant a<b\leqslant n} \ts \vp_{ab}\big(\sqrt{a(a-1)},\sqrt{b(b-1)}\ts\big);
\een
cf. \cite[Eq.~(1.3)]{n:ci}.
By Theorem~\ref{thm:fusion}, it equals $n!\ts e^{}_{\Uc}$ with
\ben
e^{}_{\Uc}=\frac{\sqrt{2}+ x_2}{2\tss\sqrt{2}}\ts\prod_{k=3}^n
\frac{x_k\big(x_k+\sqrt{k(k-1)}\ts\big)}{2\tss k(k-1)}.
\een
\eex

\section*{Declarations}

\subsection*{Competing interests}
The authors have no competing interests to declare that are relevant to the content of this article.

\subsection*{Acknowledgements}
We are grateful to Andrew Mathas for illuminating discussions
of the Jucys--Murphy method.
This project was initiated during the second and third named authors' visits
to the Shenzhen International Center for Mathematics. They
are grateful to the Center for warm hospitality.
The work of A. Molev was also supported by the Australian Research Council, grant DP240101572.

\subsection*{Availability of data and materials}
No data was used for the research described in the article.


\bigskip
\bigskip

\small

\noindent
I.K.:\newline
Shenzhen International Center for Mathematics\\
Southern University of Science and Technology, China\\
{\tt kashuba@sustech.edu.cn}

\vspace{5 mm}

\noindent
A.M.:\newline
School of Mathematics and Statistics\newline
University of Sydney,
NSW 2006, Australia\newline
{\tt alexander.molev@sydney.edu.au}

\vspace{5 mm}

\noindent
V.S.:\newline
Department of Mathematics\\
University of California at Berkeley,
Berkeley, CA 94720, USA\\
{\tt serganov@math.berkeley.edu}

\end{document}